\documentclass{amsart}
\usepackage{amssymb}
\newtheorem{theorem}{Theorem}[section]

\newtheorem{conjecture}[theorem]{Conjecture}
\newtheorem{observation}[theorem]{Observation}

\theoremstyle{definition}
\newtheorem{definition}[theorem]{Definition}
\newtheorem{example}[theorem]{Example}

\theoremstyle{remark}
\newtheorem{remark}[theorem]{Remark}
\newtheorem{question}[theorem]{Question}

\newcommand{\A}{{\mathbf A}}
\newcommand{\Q}{{\mathbf Q}}
\newcommand{\Z}{{\mathbf Z}}
\newcommand{\R}{{\mathbf R}}
\newcommand{\F}{{\mathbf F}}
\newcommand{\N}{{\mathbf Z}_{\geq 0}}

\numberwithin{equation}{section}

\begin{document}


\title[Remarks on Mazur's conjectures]{Topology of Diophantine Sets: \\
Remarks on Mazur's Conjectures}

\author[G.~Cornelissen]{Gunther Cornelissen}
\address{Ghent University \\ Department of Pure Mathematics and Computer
Algebra\\ Galglaan 2, B-9000 Gent}
\curraddr{Max-Planck-Institut f{\"u}r Mathematik \\ Vivatsgasse 7 \\
D-53111 Bonn}
\email{gc@cage.rug.ac.be}
\thanks{The first author is Post-doctoral fellow of the Fund for Scientific
Research - Flanders (FWO)}

\author[K.~Zahidi]{Karim Zahidi}
\address{Ghent University \\ Department of Applied Mathematics and
Computer Science \\ Krijgslaan 281 \\ B-9000 Gent}
\curraddr{}
\email{karim.zahidi@rug.ac.be}


\begin{abstract}
We show that Mazur's conjecture on the real topology of rational points on
varieties implies that there is no \emph{diophantine model} of the rational
integers ${\mathbf Z}$ in the rational numbers ${\mathbf Q}$, i.e., there
is no diophantine set $D$ in some cartesian power ${\mathbf Q}^i$ such that
there exist two binary relations $S,P$ on $D$ whose graphs are diophantine
in ${\mathbf Q}^{3i}$ (via the inclusion $D^{3} \subset {\mathbf Q}^{3i}$),
and such that for two specific elements $d_0,d_1 \in D$  the structure
$(D,S,P,d_0,d_1)$ is a model for integer arithmetic
$({\mathbf Z},+,\cdot,0,1)$.   

Using a construction of Pheidas, we give a counterexample to the analogue
of Mazur's conjecture over
a global function field, and prove that there is a diophantine model of the
polynomial ring over a finite field in the 
ring of rational functions over a finite field.
\end{abstract}

\maketitle

\section{Introduction}
\noindent One of the main themes in model theory is to
understand the
structure of definable sets: given a first-order language $L$
and an $L$-structure $M$, describe the 
$L$-definable subsets of $M^n$ for various $n \in \Z_{>0}$. Here, a set $S \subset M^n$ is called 
$L$-definable if there exists an 
$L$-formula $\psi ({\mathbf x})$
 with free variables ${\mathbf x}=(x_1,...,x_n)$ such that for any 
${\mathbf a} \in M^n$, 
${\mathbf a}\in S \iff M \models \psi ({\mathbf a})$.
A set
is called existentially definable (respectively, positive existentially, or diophantine) if $\psi ({\mathbf x}) $ can be taken to 
be $\exists {\mathbf b} \, \phi ({\mathbf x}, {\mathbf b})$, with $\phi $
quantifier-free (respectively, quantifier- and negation-free, or atomic).

The natural geometric examples of such structures arise as in the following
definition:

\begin{definition} \label{rl} If $R$ is a commutative ring with unit, it
admits a natural interpretation for 
any first order language $L$ of the form $L_R=(+,\cdot,=,c_i)$ where $c_i$ are
primary predicates (``constants''), less in number than $|R|$. We call
$L_R$ a 
\emph{ring language}. We define $L_{\Z} = (+,\cdot,0,1)$ and $L_t
=(+,\cdot,0,1,t)$ for any $t \in R$. \end{definition}

\begin{example} (Tarski, cf.\ \cite{Hodges:97}, pp.\ 202--206)\label{tarski1} (a) An algebraically closed field $k$ admits elimination of
quantifiers 
in the language $L_\Z$. Hence any $L_\Z$-definable subset in
$k^n$ is 
 a boolean combination of sets defined by an equation. Thus,
the
definable sets for an algebraically closed field are exactly the classical sets of algebraic
geometry -- one deduces for example that the only definable subsets of $k$ are
finite or cofinite, a fact which at first sight seems not so obvious.

(b) The field of real numbers $\R$ 
admits elimination of quantifiers
 in the language $L_{\geq}=(0,1, +, \cdot, \geq)$ of ordered fields. 
Hence every definable set in $\R^n$ is a boolean 
combination of semi-algebraic sets (i.e., solution sets to systems of 
equations of the form $f({\bf x})=0 \wedge g({\bf
x})\geq 0$).
This gives a nice description of the definable subsets of $\R$: 
they
are finite unions
 of intervals.

(c) More examples in the same vein exist, e.g., a description of
definable sets over $p$-adic fields (\cite{Mac:76}, \cite{Den:86}), 
or generalization of $(\R ,L_\geq)$ via o-minimal expansions.

(d) To give an example with a different language, existentially 
definable sets of $\Z$ in the language $(0,1,+,|)$ are unions of 
arithmetic progressions (a result of Lipshitz \cite{Lipshitz:81}). 

\end{example}

The moral is that if the (existentially) definable sets for such 
$M$ have a sufficiently easy description,  then the first-order 
(respectively, existential) theory of $M$ is decidable -- this is the case 
in the above examples. Conversely, if definable sets are combinatorially 
complicated, one expects the 
corresponding theory to be undecidable. 

\begin{example} \label{zz}
(a) Consider the rational integers $(\Z,L_\Z)$. It is impossible to describe the
$L_\Z$-definable sets of  $\Z$ 
in terms of ``classical'' sets (e.g., finite sets, arithmetic progressions, \dots). Eventually, this leads to the undecidability of the full theory of $(\Z,L_\Z)$. 

(b) The celebrated theorem of Davis, Matijasevich, Putnam and Robinson 
describes
the existentially definable sets of $\Z $: they are exactly the 
recursively enumerable sets, whose complexity outranges by far that of 
decidable (hence, certainly, of computable) sets -- and the undecidability 
of
the existential theory of $(\Z,L_\Z)$ follows. 
\end{example}

This maxim, the interplay between (un)decidability and definable sets, applies 
in particular to the field $(\Q,L_\Z)$ of rational numbers. 
The field structure of $\Q $ admits the same kind 
of ``wild'' definable sets as the integers; this follows from J.\ Robinson's theorem that $\Z $ is a definable subset of $\Q $ (\cite{Robinson:49}, theorem 3.1). The question whether the same can happen if
we restrict to the existentially definable sets is still open. 

In the next paragraph, we will present a conjecture by Mazur, which -- 
although it does not characterize the existentially definable sets of 
$\Q $ --  poses
severe 
restrictions on their real topological structure. In the subsequent 
section, we prove that this conjecture implies there is no 
``diophantine model'' (cf.\ infra) of $(\Z,L_\Z)$ in $(\Q,L_\Z)$ -- 
this generalizes Mazur's observation
that his conjecture implies that $\Z$ is not an $L_\Z$-diophantine 
subset of $\Q$. In particular, any proof of the diophantine 
undecidability of $\Q$ ``along traditional lines'' fails if 
Mazur's conjecture is true. In the final paragraph, we comment upon a 
non-archimedean version of this conjecture. 
Though most of these observations are folklore, they do not seen to have 
been written down previously.

\section{Mazur's conjectures}

In \cite{Mazur:92}, \cite{Mazur:95} and \cite{Mazur:98}, Barry Mazur has 
proposed and discussed several conjectures and questions about the 
behaviour of the set of $\Q$-rational points of a variety over 
$\Q$ under taking topological closure w.r.t.\ some metric induced by a 
valuation on $\Q$. The conjecture that we will concentrate upon 
(the weakest) is the following:

\begin{conjecture} \label{con1} \emph{(Mazur \cite{Mazur:98}, Conjecture 3)}
For any variety $V$ over $\Q$, the \textup{(}real\textup{)} topological 
closure of $V(\Q)$ in $V(\R)$ has only a finite number of real topological components.
\end{conjecture}

There is some evidence for this conjecture, especially for such $V$ which
possess special geometric properties (mostly related to the canonical 
class of $V$) -- and no counterexample to it is known. Also observe that, 
with $\Q$ 
replaced by $\R$ in \ref{con1}, the ``conjecture'' says that a real variety 
has
only finitely many real connected components. This holds true; it could be deduced
from Tarski's results -- there is 
even an explicit
bound on the betti numbers of $V(\R)$, the so-called Milnor-Thom theorem
(cf. \cite{Milnor:64}).

\begin{example} (a) Conjecture \ref{con1} is true for curves $V$. One can 
assume
$V$ to be projective and non-singular. 
The case where $V$ has genus $ g \geq 2$ is settled by Faltings's theorem, 
which says that $V(\Q)$ is a finite set (\cite{Faltings:83}). 
If $V$ has genus 0, then either $V(\Q )$ is empty, or $V$ is 
$\Q$-birational to $\A^1$, and $\A^1(\Q)$  is topologically dense 
in $\A^1(\R)$. Finally, assume that $V$ has genus 1. It is known that $V(\R)$ is
isomorphic
to the ``circle group'' $\R / \Z$ or to $\R / \Z \times \Z / 2$ (see
\cite{Sil:95}, V). Every proper 
closed subgroup of the circle group is
finite (see \cite{Hoc:88}, theorem 1.34). Hence, if $V(\Q)$ is not finite, 
then
it is dense in every component of $V(\R)$ that it intersects.           

(b) To provide a higher dimensional example, let $V$ be a variety 
satisfying weak approximation (i.e., such that $V(\Q) \hookrightarrow 
\prod V(\Q_p)$ is dense). Then the conjecture holds true for $V$. 
This holds, e.g., if $V$ is a smooth complete intersection of two 
quadrics in projective space of dimension at least 5 (cf.\ \cite{Mazur:92}). 
\end{example}

\begin{remark}
Mazur has made even stronger conjectures, some of which had to be slightly 
modified, due to the construction of a counterexample by 
Colliot-Th{\'e}l{\`e}ne, Skorobogatov and Swinnerton-Dyer 
(\cite{Colliot:97}). For an extensive (unsurpassable) exposition and more 
examples, we refer to the original sources \cite{Mazur:92}, \cite{Mazur:95} 
and \cite{Mazur:98}. 
We will concentrate on the model-theoretical aspects of the conjectures, 
which are already present in \ref{con1} -- but let the reader be warned 
about making 
too bold generalizations of \ref{con1}. A non-archimedean version will be
considered in the last paragraph of this paper. 
\end{remark}

\begin{remark}
The $(\Q,L_\Z)$-existentially definable subsets, in the sense of the 
introduction, are precisely images of projections from $V(\Q)$ to affine
space $\A^n_\Q$ for various $V$ and $n$.
 
A more model-theoretic conjecture would be that \emph{the real topological 
closure of a $(\Q,L_\Z)$-existentially definable set 
is an $(\R,L_{\geq})$-definable set} (i.e., a semi-algebraic set). This 
implies \ref{con1}, since 
a semi-algebraic set has only finitely many components. 

 We do not know whether conjecture \ref{con1} is equivalent to this 
statement. Note that J.\ Robinson's argument (in \cite{Robinson:49}) shows that it 
is wrong when the word ``existentially'' is erased.  
\end{remark}

\section{Diophantine models of $\Z$ in $\Q$}

Mazur has observed that conjecture \ref{con1} implies that $\Z$ is not 
diophantine in $\Q$ in the language $L_{\Z}$; indeed, if $\Z$ would arise as the projection of $V(\Q)$ for
some variety $V$, then since $\Z$ has infinitely many real components and 
the projection is continuous, the same would hold for $V(\R)$. 

However, many proofs of the undecidability of the diophantine theory of
structures $(R,L_R)$ as in (\ref{rl}) do not give that $\Z$ is a diophantine
\emph{subset} of $R$, but rather produce a \emph{diophantine model} of
$(\Z,L_\Z)$ in $(R,L_R)$, in the sense of the following 
definition:

\begin{definition} 
A model $(M,L,\phi)$ is a triple consisting of a first order language $L$ 
which consists of $i$-ary predicates $\{ P_{i,\alpha} \}$, a set $M$ and an interpretation
$\phi$ of $L$ in $M$ (we will often leave out $\phi$ of the notation). Note
that any cartesian power $M^k, \ (k \geq 1)$ is likewise a model for $L$
via ``diagonal interpretation''.

We say that a model $(M',L'=\{P'_{i,\alpha} \},\phi')$ admits a
\emph{diophantine model} in $(M,L,\phi)$ if there exists a set-theoretical
bijection between $M'$ and a subset of some cartesian power $M^k \ (k \geq
1)$, such that the image is diophantine, and such that the induced
inclusions $\phi'(P'_{i,\alpha}) \subseteq M^{ik}$ are diophantine.

 A similar notion of
\emph{\textup{(}positive\textup{)} existential model} exists.
\end{definition}

\begin{example} \label{ex}
(a) If $(M^2,L)$ admits a diophantine model in $(M,L)$, then the latter
structure is said to admit \emph{diophantine storing} (cf.\
\cite{Cornelissen:99}). This is true, for example, for $(\Z,L_\Z)$. For
non-algebraically closed rings $(R,L_R)$ admitting diophantine storing, one
can always choose $k=1$ in the above definition. For if $(M',L',\phi ')$ admits
 a diophantine model in $(R^2, L_R)$ and $(R,L_R)$ admits diophantine storing,
 then $(M',L',\phi ')$ admits a diophantine model in $(R,L_R)$ (since 
conjuntions of
 diophantine formulas are again diophantine if the quotient field of
 $R$ is not algebraically closed -- cf. \cite{Cornelissen:99}, \S 3).

(b) Typically, diophantine models of the integers $(\Z,L_\Z)$ in ring languages $(R,L_R)$ arise in the following way: a commutative algebraic group $G$ (e.g., the multiplicative group of a quadratic ring, or an elliptic curve) is assumed to have rank one 
over $R$, and the set $\Z$ has a diophantine model as the
$R$-rational points $G(R)$ of $G$ - the relation ``addition'' is 
automatically mapped to a diophantine subset of $G^3(R)$, since the group 
law on $G$ is a morphism.   The most problematic point is defining the 
relation ``multiplication''. For an example, consider the proof that 
$(\Z,L_\Z)$  admits a diophantine model in $(R:=S[t],L_t)$ for any
commutative unitary domain $S$ of characteristic zero, see Denef 
\cite{Denef:78}. 
He
takes for $G$ the torus ${\bf G}_{m,R[\sqrt{\Delta}]}$ of discriminant 
$\Delta=t^2-1$, which is non-split over $R$; $G(R)$ has rank one: any 
$R$-point is given by a
solution $(x_n,y_n)$ to the Pell-equation $X^2-\Delta Y^2=1$ (i.e., a 
power $u^n=x_n + y_n \sqrt{\Delta}$ of the fundamental unit 
$u=t+\sqrt{\Delta}$). Multiplication $(x_r,y_r) \cdot (x_s,y_s) = (x_n,y_n)$ 
is defined by saying that $f:=y_n-y_r \cdot y_s$ has a zero at 
$t=1$, i.e., $(\exists h \in R)(f=(t-1)h).$ 

(c) It is not known whether, if the ring $R$ contains $\Z$, the set $\Z$ itself
is a diophantine subset of $R$ whenever $(\Z,L_\Z)$ admits a diophantine model
in $(R,L_\Z)$. 
\end{example}

The following result formalizes the technique of proof of many
undecidability results:

\begin{observation}
Assume that
$R$ is as in \ref{rl}, and there is a polynomial whose coefficients belong
to $\phi(L)$ but that has no zero in the fraction field of $R$. If
$(M',L')$ has an undecidable diophantine theory, and admits a diophantine
model in $(R,L_R)$, then the diophantine theory of $(R,L_R)$ is
undecidable. \qed
\end{observation}

\begin{remark}
Without any restrictions on $R$, if $(M',L')$ has an undecidable (positive)
existential theory and admits a (positive) existential model in $(R,L_R)$,
then
the (positive) existential theory of $(R,L_R)$ is undecidable.
\end{remark}

The technique of many undecidability proofs for rings $(R,L_R)$ is based on
the fact that, via a construction as in (\ref{ex}(b)), one can find a diophantine model of the integers $(\Z,L_\Z)$ in
$(R,L_R)$, and then rely on the fact that the diophantine theory of the
integers is undecidable (\cite{Mat:70}, \cite{Davis:76}).  
It has been suggested that, with this more flexible definition, one would
be able to find a diophantine model of the integers in the rationals:

\begin{question} \label{ZinQ}
Does $(\Z,L_\Z)$ admit a diophantine model in $(\Q,L_\Z)$?
\end{question} 

However, even this is impossible if we assume Mazur's conjecture:

\begin{theorem} \label{th1}
Mazur's conjecture \ref{con1} implies that there is no diophantine model 
of $(\Z,L_\Z)$ in $(\Q,L_\Z)$. 
\end{theorem}

\begin{proof}
Assume that there is such a diophantine model $(D,L_D)$, with $D \subseteq
{\bf Q}^k$. Then there is an affine variety $V$ over $\Q$  admitting
a finite morphism  $f:V_{\Q} \rightarrow {\mathbf A}^k_\Q$ defined over
$\Q$ such
that $f(V(\Q))=D$.

If $D$ is discrete (i.e., infinite and totally disconnected in the real topology),
the traditional proof applies: the real topological closure of $V(\Q)$ in
$V(\R)$ is also mapped to $D$ by
$f$, and hence it has infinitely many components.  

If $D$ is not discrete (which seems to be the case for the typical infinite 
diophantine set in $\Q$, say, the set of squares), then we show that one can 
select (in a computable way) a discrete subset $\tilde{D}$ of $D$. Then 
the above proof, applied
 to $\tilde{D}$, gives the result.

Here are the details of the construction of $\tilde{D}$. We only have to
treat the case where the real topological closure $\bar{V}$ of $V(\Q)$ has
only finitely many connected components. Since $f$ is continuous, the mean
value theorem implies that $f(\bar{V})$ is the union of finitely many
closed subsets in $\R^k$. In particular, the topological closure $\bar{D}$
of $D$ contains finitely many closed subsets, and since $D$ is infinite,
one of these subsets, say, $D_0$, is not a point. By composing $f$ with a
suitable $\Q$-rational projection $\pi \ : \  {\mathbf A}^k_\Q \rightarrow
{\mathbf A}^1_\Q$ which does not map $D_0$ to a point, we may assume $k=1$.
By composing with a fractional linear transformation defined over $\Q$, we
may assume $\pi(D_0)$ to be the unit interval $ I = [0,1]$.  Let $d_n$ be
the element of $d$ corresponding to $n \in \Z$.
Let us consider the set
$$ \tilde{Z} = \{ n \in \Z \ | \ \frac{1}{2j+1} \leq \pi(d_n) \leq
\frac{1}{2j} \mbox{ for some } j \in \Z_{>0} \}. $$ 
The set $\tilde{Z}$ is Turing computable (since $D=\{ \pi(d_n) \}$ is a
listable subset of $\Q$, it is easy to write a Turing program to check the
inequalities), hence it is recursively enumerable (by Kleene's normal form
theorem, cf. \cite{Soare:87}, 2.3-2.4), so by \cite{Davis:76}, it is
diophantine in $(\Z,L_\Z)$.
Also, $\tilde{Z}$ is infinite, since $\pi(D) \cap I$ is dense in $I$. 
 We now set 
$$\tilde{D} = \{ d_n \ | \ n \in \tilde{Z} \}. $$
By construction, the set $\tilde{D}$ is diophantine in $(D,L_D)$, and hence
a fortiori in $(\Q,L_\Z)$. So there exists a variety $\tilde{V}$ over $\Q$ and
a $\Q$-morphism $\tilde{f}: \tilde{V} \rightarrow  {\mathbf A}^1_\Q$ such
that $\tilde{f}(\tilde{V}(\Q)) = \tilde{D}$. However, the real closure of
the set $\tilde{D}$ has infinitely many connected components in the real
topology by construction. Hence the same
holds for $\tilde{V}(\Q)$, contradicting Mazur's conjecture.  
\end{proof}

\section{Non-archimedean aspects of Mazur's conjectures} 

In (\cite{Mazur:98}, II.2), Mazur has devised a conjecture of the above
type which applies to any completion of a number field, not just an
archimedean one. 
As it makes sense for any global field, we formulate it as follows:

\begin{question}  \label{q1} Let $V$ be a variety over a global field $K$,
$v$ a valuation on $K$, and $K_v$ the completion of $K$ w.r.t.\ $v$. For
every point $p \in V(K_v)$, let $W(p)$ be the Zariski closure of $\bigcup
(V(K) \cap U)$, where $U$ ranges over all $v$-open neighbourhoods of $p$ in
$V(K_v)$. Is the set $\{ W(p) \ : \ p \in V(K_v) \}$ finite?
\end{question}

In our next theorem, we observe that the answer to this question is
negative in positive characteristic:

\begin{theorem} \label{wrong}
Let $K=\F_q(t)$ be the rational function field over a finite field $\F_q$
of positive characteristic $p>0$, and
let $v$ be the valuation corresponding to the place $t^{-1}$ of $K$
\textup{(}i.e., $v(a)=q^{\deg(a)}$ for $a \in \F_q[t]$\textup{)}. Then
there is a variety $V$ over $K$ for
which the answer to question \ref{q1} is negative. 
\end{theorem}

\begin{proof}
In \cite{Pheidas:91} (lemma 1) Pheidas proves that, for $p>2$, projection
onto the  $x$-coordinate of the $K$-rational points of the space curve
$V_p$ given by $$ V_p \ : \ \ x-t=u^p-u, \ x^{-1}-t^{-1} = v^p-v $$
gives the set 
$ D_p = \{ t^{p^s} \ | \ s \in \Z_{\geq 0} \}.$ For $p=2$, Videla
(\cite{Videla:94}) proved that the set $D_2$ is the projection onto the
$x$-coordinate of  
$$ V_2 \ : \ x+t=u^2+u, \ u=w^2+t, \ x^{-1}+t^{-1}=v^2+v, v=s^2+t^{-1}. $$
Already the sets $W(p)$ for $p \in V(K)$ are disjoint, since their
$x$-coordinates are separated ($v(t^{p^s}-t^{p^r})>1$ for
all $r \neq s$). This gives a negative answer to question \ref{q1}. 
\end{proof}

Thinking of the analogy between function fields and number fields, one can 
ask for the strict analogue of question \ref{ZinQ} for global 
function fields.The answer to it is \emph{positive}:

\begin{theorem}
For any prime power $q$, $q=p^n$, $p>0$, the polynomial ring $(\F_q[t],L_t)$ admits a
diophantine model in the ring of rational functions $(\F_q(t),L_t)$.
\end{theorem}

\begin{proof}
The proof is a bit indirect: we show that the polynomial ring has a
diophantine
model in the positive rational integers, and the latter has a diophantine model in
the field of rational functions. 

More precisely, $\F_q[t]$ is a recursive ring (cf.\ Rabin \cite{Rabin:60}),
because $\F_q$ is recursive (since finite), and hence the same holds for
the polynomial ring over $\F_q$ (cf.\ Fr{\"o}hlich and Sheperdson
\cite{Frohlich:56}). So there exists an injective map $ \theta \ : \
\F_q[t] \rightarrow \N$
such that the graphs of addition and multiplication are recursive on 
$\N$, and hence $(\theta(\F_q[t]),\theta(L_t))$ is
a diophantine model of $(\F_q[t],L_t)$ in $(\N, L_\Z)$. 

For the second step, we first recall a construction of Pheidas and Videla (\cite{Pheidas:91}, \cite{Videla:94}).
Let $v$ denote the $t$-valuation on $\F_q(t)$, i.e., $v(x)$ is the order of $x$ at zero. For any $k\in {\mathbf Z}_{\geq 0}$, let $[k]$ denote the equivalence class
of elements 
$x \in \F_q(t)$ with $v(x)=k$.  For positive integers $a$ and $b$, let 
$a\, |_p\, b$ denote the relation $(\exists n \in \N)(a=bp^n)$. 

Consider the structure $S=(\Z_{\geq 0}, (+, |_p, 0,1))$. Firstly, 
multiplication is diophantine in $S$ (\cite{Pheidas:87}, corollary on p. 529). 
Secondly, the set of equivalence classes $[k]$ as above is a model for $S$ in which the relations
in $S$ can be defined by diophantine formulas in $(\F_q(t),L_t)$ between arbitrary representatives of the 
equivalence classes in $\F_q(t)$. 
We conclude that for arbitrary elements $x,y,z \in \F_q(t)$
 the relations $[v(x)]=[v(y) + v(z)]$ and 
$[v(x)]=[v(y)\cdot v(z)]$ are diophantine in 
${\mathbf F}_q(t)$.


The problem with this encoding is that we do not know the existence of a diophantine set 
 in $\F_q(t)$ which contains exactly one representative for each such equivalence 
class. We fix this problem as follows. We know from the proof of theorem \ref{wrong} that 
the set 
$D_p=\{ t^{p^k}, \; k\in \Z _{\geq 0} \}$ is diophantine in $(\F_q(t),L_t)$, and this
 will be our model.

To define addition and multiplication on elements of this set, we introduce the
 following switching between
 $t^{p^k}$ and $[k]$: the set $\{ (k,p^k), \; k\in \Z _{\geq 0} \}$ 
is recursively enumerable in $\Z _{\geq 0}^2$, so by Matijasevich's 
theorem, it is diophantine in $\N$. Then, by the aforementioned results, 
the set ${\mathcal E}=\{ ([k],[p^k]), \; k\in \Z _{\geq 0} \}$
is diophantine over $(\F _q(t),L_t)$. 

For the function symbols $R\in \{ +, \cdot \}$ on $\N$, we let the corresponding symbol $\tilde{R}$ for $x,y,z \in D_p$ be
defined by
\begin{align*} 
z=x\tilde{R} y \Leftrightarrow & (\exists \; x_1,y_1,z_1 \in \F_q(t)) 
 (((x_1,x),(y_1,y),(z_1,z)) \in {\mathcal E}^3 
 \wedge \; [R(x_1,y_1)]=[z_1]) .\end{align*}
For $R \in \{ +, \cdot \}$, the righthand side of the equivalence is diophantine in $(\F_q(t),L_t)$ by what 
we have said before and the fact that for any two elements 
$w_1,w_2 \in \F _q(t)$, the statement $[w_1]=[w_2]$ 
is equivalent with $(v(w_1/w_2) \geq 0) \wedge (v(w_2/w_1) \geq 0)$, which is 
diophantine by \cite{Pheidas:91} 
and \cite{Videla:94}.

Finally, $(D_p,\tilde{+},\tilde{\cdot},t,t^p)$ is a diophantine model of 
$(\Z, L_{\Z })$ in 
$(\F _q(t), L_t)$. This finishes the proof of the theorem.
\end{proof}

Of course, the above theorem still does not settle the following problem:

\begin{question}
Is the polynomial ring $\F _q[t]$ a diophantine \emph{subset} of the field of 
rational functions $\F _q(t)$?
\end{question}

\bibliographystyle{amsplain}

\end{document}